\documentclass[11pt,a4paper]{article}
\usepackage{graphicx}
\usepackage{cite}
\usepackage{amssymb,amsmath}
\usepackage{layout}
\usepackage{subcaption}
\usepackage{amsthm}
\usepackage{algorithm}
\usepackage{algpseudocode}
\usepackage[makeroom]{cancel}
\usepackage[dvipsnames]{xcolor}

\newtheorem{proposition}{Proposition}[section]

\DeclareGraphicsExtensions{.eps,.jpg,.pdf,.ps}

\graphicspath{{./Figs/}}

\begin{document}

\bibliographystyle{unsrt}

\def\EQN{\eqno}
\def\beq{\begin{equation}}
\def\eeq{\end{equation}}
\def\beqn{\begin{eqnarray}}
\def\eeqn{\end{eqnarray}}
\def\bef{\begin{figure}}
\def\eef{\end{figure}}
\def\del{\partial}
\def\Psiz{\Phi}
\def\psinc{\psi_{inc}}
\def\vnc{\v_{inc}}
\def\r{{\bf r}}
\def\rr{{\bf r'}}
\def\d{{\bf d}}
\def\h{{\bf h}}
\def\A{{\bf A}}
\def\B{{B}}   
\def\bea{\begin{aligned}}
\def\eea{\end{aligned}}
\def\zz{Z}
\def\Arcsin{{\rm Arcsin}}
\def\Arctan{{\rm Arctan}}
\def\t{\theta_n}
\def\tt{\frac{\theta_n}{2}}
\def\hmax{h_{\max}}
\def\hmin{h_{\min}}
\def\P{{P}}   
\def\Q{{Q}}     
\def\Acal{{\cal{A}}}
\def\Ecal{{\cal{E}}}
\def\Fcal{{\cal{F}}}
\def\Gcal{{\cal{G}}}
\def\Scal{{\cal{S}}}
\def\Hcal{{\cal{H}}}
\def\Bcal{{\cal{B}}}
\def\Pcal{{\cal{P}}}
\def\Qcal{{\cal{Q}}}
\def\Ical{{\cal{I}}}
\def\Jcal{{\cal{J}}}
\def\Lcal{{\cal{L}}}
\def\Rcal{{\cal{R}}}
\def\Mcal{{M}}  
\def\Ncal{{\cal{N}}}
\def\Einc{{\u}}
\def\Escat{{E_s}}
\def\Epartial{\frac{\partial E}{ \partial n}}

\def\Sn{{S_n}}
\def\Sone{{S_{n+1}}}
\def\Stwo{{S_{n+2}}}
\def\v{{\bf v}}
\def\u{{\bf u}}
\def\w{{\bf w}}

\def\rbar{{\bf r}}
\def\lbar{{\bar{\lambda}}}
\def\J{{\bf J}}
\def\t{{\bf t}}

\parindent 0 true cm

\title{Invariant subspaces and nest algebras\\
in left--right operator splitting}
\author{Mark Spivack\\[4pt]
{\small Department of Applied Mathematics and Theoretical Physics,}\\
{\small University of Cambridge, UK}\\
{\small\texttt{ms100@cam.ac.uk}}}

\date{\today}
\maketitle

\begin{abstract}
The left-right (L-R) splitting method for boundary integral equations in rough-surface scattering writes the integral operator as the sum of its lower- and strictly upper-triangular parts. This note examines the operator algebra formed by the lower-triangular ($L$-type) parts. Their weak closure is a nest algebra; second-kind $L$-type operators with Hilbert-Schmidt kernels have $L$-type inverses; and Arveson's distance formula expresses the error of the best possible forward-only approximation in terms of the backward coupling across a single cut. Bounds and numerical results for the 2D double-layer kernel suggest that truncation is close to optimal, and that the worst-case input field is a predominantly left-going grazing wave rather than backscatter of a right-going incident field. For a hard surface the inverse $L$-operator is given, up to a second-order correction, by four times the $L$-operator of the mirror-image surface, and several open questions are posed.	

\end{abstract}

\section{Introduction}
\label{sec:intro}

This note explores some operator-algebraic concepts associated with
integral operators and invariant subspaces arising in operator splitting
for surface scattering, applies the resulting analytical tools to the method,
and poses further open questions.

{\bf Motivation:} Suppose $A$ is an integral operator on a
Hilbert space $H$ of functions on $\mathbb{R}^n$ or on a compact
subset, and that we wish to solve $A\u = \v$ to find $\u ~=~ A^{-1}\v$, where
$\v$ is a known function (`incident wave') and $\u$ is unknown (`surface field').
Where this solution exists it can formally be written as a series:
\beq
A^{-1}\v ~=~ L^{-1} \sum_{n=0}^\infty B^n \v
\label{(3.6new)}
\eeq
which may be truncated, where $\B = - R L^{-1}$. This is the
left--right (L-R) splitting method 
\cite{spivack2001validation,spivack2004left}.
The series (\ref{(3.6new)}) is analogous to the Bremmer coupling
series for layered 
media \cite{bremmer1951wkb,dehoop1996bremmer}, and closely related to
the method of ordered 
multiple interactions and the forward--backward method for rough
surfaces \cite{kapp1,tran1997calculation,pino}. 
The operators $L$ and $R$ are respectively
(borrowing terminology from the finite-dimensional case)
the lower- and strictly upper-triangular parts of $A$; they are
defined more precisely in 
Section~\ref{sec:setting}.
We call operators of this lower-triangular form \emph{$L$-type}: integral
operators, possibly with an added multiplication operator, whose
kernels vanish for 
$x'>x$, so that $(L\psi)(x)$ depends only on the values of $\psi$ to
the left of $x$. 
It is convenient to think of any point $a$ as a {\sl cut}, splitting
the surface into 
the part $x\le a$ and the part $x>a$. An $L$-type operator then
carries information 
forward across every cut, but none backward.
The notation $\Lcal$ for the set of such operators is also introduced in
Section~\ref{sec:setting}.

Our motivation arises partly
from viewing the $L$-operators as a directional flow on the
surface, using this to characterise the inverses $L^{-1}$
and to bound the errors in forward-only approximations.

\medbreak
{\bf Summary of results.} Our main observations are as follows.
\begin{itemize}
\item All $L$-type operators leave invariant a continuous nest $\Ncal$ of
subspaces, consisting of the functions that vanish to the left of each point, and
their weak closure is the {\sl nest algebra} $\mathrm{Alg}\,\Ncal$ of all bounded
operators leaving $\Ncal$ invariant; moreover $\Lcal$ is triangular in the sense of
Kadison and Singer, and $\mathrm{Alg}\,\Ncal$ is maximal triangular (Section~\ref{sec:nest}).
\item The inverse of an invertible $L$-type operator need not in
  general be of $L$-type itself. However, for second-kind $L$-operators
  whose multiplication term is invertible 
in $L^\infty$ (bounded away from zero) and whose kernel is Hilbert--Schmidt (for instance on a bounded interval), the inverse is again of $L$-type; this supports forward marching (Section~\ref{sec:nest}).
\item Arveson's distance formula gives the error of the best possible forward-only
approximation to $A$, valid for every splitting and not only
triangular truncation, in terms of
the backward coupling across a single cut. For the double-layer kernel of a rough
surface we bound this quantity and compute it numerically. In our examples plain
truncation is close to optimal, and the worst case on all of $L^2$ is a predominantly
left-going wave rather than backscatter of the incident field (Section~\ref{sec:arveson}).
\item
For a hard surface $h$, the inverse $L$-operator is closely related to the forward $L$-operator for the
mirror surface $z=-h(x)$.
Specifically, $L_h^{-1}$ is exactly $4L_{-h}$, times a correction factor that differs from
the identity only at second order in the kernel. For gently sloping surfaces, the
mirror-surface operator alone therefore gives the inverse to first order (Section~\ref{sec:mirror}).
\end{itemize}
The operator-algebraic theorems used here are largely classical;
the purpose of the note is to identify their implications for left--right
splitting of boundary-integral scattering equations. In particular, we use
Arveson's distance formula to assess the limitations of forward-only
approximations, examine the directionality of the fields that produce
the strongest backward coupling across a cut, and
derive a mirror-surface identity for the
Neumann forward operator.
The numerical observations
suggest possible directions for further developments and open questions
which are described in Section~\ref{sec:questions}.

\section{Setting and notation}
\label{sec:setting}

We take the wave functions in our problem to lie in a Hilbert space $H$, usually
$L^2(\Omega)$, where $\Omega$ is $\mathbb{R}$ or $\mathbb{R}^2$ or a compact subset
of these. For weakly singular or first-kind operators the Sobolev spaces
$H^k(\Omega)=W^{k,2}(\Omega)$ are often more natural (see the Remark on boundedness
below).
Denote by $B(H)$ the set of bounded linear operators on $H$. We consider wavefields defined on surfaces expressible as
functions of horizontal coordinates $x$ or $(x,y)$. For
now we restrict attention to $x$, i.e.\ the 2D problem.

\medbreak
\textbf{First- and second-kind integral equations.}
An integral equation of the first kind has the form $Ku=f$, where $K$
is an integral operator. An equation of the second kind has the form
$(M_c+K)u=f$, where $M_c$ denotes multiplication by a coefficient $c$
bounded away from zero, i.e.\ invertible in $L^\infty$ (often a constant, such as $\frac12$ in boundary-integral
equations). First-kind operators may be smoothing; when $K$ is compact, as
it is for weakly singular kernels such as the logarithmic one, its inverse cannot
be bounded on $L^2$, even when the equation between
suitable Sobolev spaces is well posed. The presence of a
multiplication term does 
not by itself guarantee invertibility in the second-kind case. For the
$L$-operators considered below, the Volterra structure follows automatically
from the support condition $x'<x$; in Proposition~\ref{prop:inverse} it is
the regularity assumptions that ensure invertibility.

\medbreak
{\bf Integral operators:} We choose a kernel class $\mathcal{K}$ on
$\Omega\times\Omega$ 
whose elements define bounded integral operators on $H$ and which is
closed under 
kernel composition and adjoints; for instance, Hilbert--Schmidt kernels
($K\in L^2(\Omega\times\Omega)$) or Schur-class kernels (defined
below). Consider the 
set $\Jcal\subset B(H)$ of integral operators of the form
\beq
\Jcal = \{ A_K: A_K\psi(x) \equiv \int \psi(x')K(x,x')dx' +\alpha(x)\psi(x), ~~~
~~~\psi \in H \},
\eeq
for some kernel $K\in\mathcal{K}$ and $\alpha\in L^\infty(\Omega)$.
The term $\alpha$ is equivalent to allowing delta functions in the
kernel, i.e.\ kernels may contain $\alpha(x)\delta(x-x')$,
and thus $\Jcal$ contains the identity and all multiplication operators.
This is also physically natural, allowing for variable diagonal terms.
With this choice of $\mathcal{K}$, $\Jcal$ is a unital normed operator
algebra, with multiplication given by composition of mappings (though
it is not in 
general norm-closed in $B(H)$).
Explicitly, if $J, K \in \mathcal{K}$
are two kernels, their composition $L$ is a kernel where
\beq
L(x,x'') = \int K(x,x')J(x',x'')\,dx'
\eeq
so
\beq
\begin{split}
A_K A_J \psi(x) &= \int \left\{\int \psi(x'')J(x',x'')\,dx'' \right\}
K(x,x')\,dx' \\
&= \int  \psi(x'')\left\{ \int J(x',x'') K(x,x')\,dx'\right\} dx''
~=~ A_L\psi(x) .
\end{split}
\label{kernelproduct}
\eeq

\medbreak
{\bf Left-type operators:} Now consider the mapping $A\rightarrow L$ where, for an
integral operator $A\psi(x) = \int \psi(x')K(x,x')\,dx'$ with $K\in\mathcal{K}$,
\beq
L\psi(x) = \int_{x'< x}
\psi(x')K(x,x')\,dx' ~=~
\int \psi(x')\Theta(x-x')K(x,x')\,dx' ,
\eeq
where $\Theta$ is the Heaviside step function (for conditionally convergent or
principal-value integrals, the integral is interpreted as
$\lim_{\epsilon\rightarrow 0^+} \int_{x'< x-\epsilon}\psi(x')K(x,x')\,dx'$). Any
diagonal term $\alpha(x)\psi(x)$ of $A$ is kept in $L$.
We denote this mapping by $\Pcal(A)=L$, so that $\Pcal:\Jcal \rightarrow \Jcal$.
We call operators of the form $L=\Pcal(A)$, $A\in\Jcal$, {\sl $L$-type} (or
{\sl left-type}) operators, and write $\Lcal=\Pcal(\Jcal)$ for the set of such operators.
Their kernels are supported in $x'\le x$, so the value of $L\psi$ at $x$ depends only on
$\psi$ to the left of $x$.
The set $\Lcal$ is itself an algebra: sums and products of $L$-type operators are
again $L$-type, since if two kernels vanish for $x'>x$ then so does their
composition (\ref{kernelproduct}), while the multiplication terms simply multiply.

The mapping $\Pcal$ is a projection, since $\Pcal^2 = \Pcal$:
\beq
\begin{split}
\Pcal^2(A)\psi ~=~ \Pcal(L)\psi &= \int
\psi(x')\Theta^2(x-x')K(x,x')\,dx' \\
&= \int
\psi(x')\Theta(x-x')K(x,x')\,dx' ~=~ \Pcal(A)\psi .
\end{split}
\eeq
Furthermore, the operator $\Qcal = 1-\Pcal$ is a projection, since
$\Qcal^2=\Qcal$ (immediate),
and $\Qcal(A)=R$. We call $\Qcal$ the complement of $\Pcal$, and by
slight abuse of notation we may denote it $\Pcal^\perp$. Operators of the form
$R=\Qcal(A)$ are {\sl $R$-type}, and we write $\Rcal=\Qcal(\Jcal)$.

{\bf Norms of $\Pcal$ and $\Qcal$.} Since $\Pcal$ acts on an operator by
multiplying its kernel by $\Theta(x-x')$, we use the same symbols $\Pcal$ and $\Qcal$
for the induced maps on kernels, $K\mapsto\Theta(x-x')K$ and
$K\mapsto\big(1-\Theta(x-x')\big)K$. Both $\Pcal$ and $\Qcal$ have norm 1
in two useful norms:
\begin{itemize}
\item On Hilbert--Schmidt kernels,
$\Pcal$ is the orthogonal projection of $L^2(\Omega\times\Omega)$ onto kernels supported
in $x'<x$, and $\Qcal$ is the complementary orthogonal projection.
\item On the Schur class
($\max\{\sup_x\int|K|\,dx',\ \sup_{x'}\int|K|\,dx\}<\infty$), $\Pcal$
  and $\Qcal$ are 
contractive since $|\Pcal(K)|\le|K|$ and $|\Qcal(K)|\le|K|$
pointwise. The Schur class 
is a Banach algebra acting on every
$L^p$, and contains weakly singular Green's-function kernels on
bounded intervals. 
\end{itemize}
In either of these norms
\beq
||\Pcal(A)|| \le ||A||, ~~~~~~~ ||\Qcal(A)|| \le ||A||.
\label{contractive}
\eeq
In the Hilbert--Schmidt norm the inequalities are in fact strict unless
$\Pcal(A)=A$ or $\Qcal(A)=A$ respectively, since
$\|A\|_{HS}^2=\|\Pcal(A)\|_{HS}^2+\|\Qcal(A)\|_{HS}^2$. 
In the Schur norm, which is a supremum,
the inequalities need not be strict.
The norm-one bounds (\ref{contractive}) imply that
$\mathrm{ran}(\Pcal)$ and $\ker(\Pcal)$ are mutually Birkhoff--James
orthogonal (that is, $\Vert{}L\Vert{} \le \Vert{}L+\lambda R\Vert{}$
and $\Vert{}R\Vert{} \le \Vert{}R+\lambda L\Vert{}$ for all
$L\in\mathrm{ran}(\Pcal)$, $R\in\ker(\Pcal)$, and
$\lambda\in\mathbb{C}$), which partly motivates the notation
$\Pcal^\perp$ (see e.g.\ \cite{sain2018complete}).
The bounds (\ref{contractive}) do not hold in the operator norm on $B(H)$:
triangular truncation is unbounded in operator norm
(for $n\times n$ matrices its norm grows like $\log n$
\cite{kwapien1970main}; for convolution kernels it acts like a Hilbert
transform on the symbol).

{\sl Remark (boundedness).} We wish to allow kernels involving Green's functions,
which are not in general essentially bounded and need not give rise to bounded
operators on $L^2$. Integral operators with weak singularities map $W^{k,p}$ into
$W^{k+1,p}$ boundedly \cite{witte1999integral}, which may be sufficient for us. If
we extend to unbounded operators and/or non-closed inner product spaces we may still
be able to discuss nest algebras and their properties in a generalised sense; for
first-kind equations such an extension is unavoidable (Section~\ref{sec:nest}).

\section{Nest algebra, inverses, and triangularity}
\label{sec:nest}

{\bf The subspace nest.}
Consider the set $\Ncal$ of closed subspaces $\{ \Mcal_a: a \in
\mathbb{R}\cup\{\pm\infty\}\}$ of
$H=L^2(\mathbb{R})$ consisting of functions which vanish for $x \le a$, i.e.\
$\Mcal_a =  \{\psi:  \psi(x)=0 ~~\mbox{for}~~ x \le a  \} $, together with
$\Mcal_{-\infty}=H$ and $\Mcal_{+\infty}=\{0\}$ (on a bounded interval $[0,\ell]$,
$a$ ranges over $[0,\ell]$ with $\Mcal_0=H$ and $\Mcal_\ell=\{0\}$).
Each $\Mcal_a$ is a subspace, and $\Mcal_b \subseteq \Mcal_a$
for $a \leq b$, so
that these subspaces are nested (totally ordered by inclusion).
In fact $\Ncal=\{ \Mcal_a\}$ is a {\sl complete lattice} (every subset has
a supremum and an infimum in $\Ncal$),
so $\Ncal$ is a \textit{subspace nest}. Similarly we can consider
the orthogonal complements $\Mcal_a^\perp$, which are of the
form $\Mcal_a^\perp =  \{\psi:  \psi(x)=0 ~~\mbox{for}~~ x > a  \} $, and
denote the set $\{\Mcal_a^\perp\}$ by $\Ncal^\perp$.

Let $P_a$ be the orthogonal
projection onto $\Mcal_a$, so that $P_a^\perp=I-P_a$ projects onto the functions
supported in $x\le a$.
For any $\psi \in \Mcal_a$ and integral operator $A$ with
lower-triangular part
$L=\Pcal(A)$, it is apparent
that $L\psi$ is also in $\Mcal_a$. Therefore $\Mcal_a$ is invariant under
any such $L$; in other words $L P_a = P_a L P_a$ for all $a$ and all $L \in \Lcal$.
An analogous result applies to the upper-triangular part
$R=\Qcal(A)$ with respect to $\Ncal^\perp$; that is,
$R P_a^\perp = P_a^\perp R P_a^\perp$ for all $a$ and all $R\in\Rcal$.

{\bf Adjoints.} If $A$ is defined by
$A\psi = \int \psi(x')K(x,x')dx'$, then the adjoint $A^*$ is given by
\beq
A^*\psi(x) = \int \psi(x') \overline{K(x',x)} dx' .
\eeq
Writing $L$ in terms of the Heaviside function,
\beq (L\psi)(x) =  \int\psi(x')~\Theta(x-x')~K(x,x')~dx' ,
\eeq
the adjoint is
\beq (L^*\psi)(x) =  \int\psi(x')~\Theta(x'-x) ~\overline{K(x',x)}~dx'
~=~ \int_x^\infty \psi(x') ~\overline{K(x',x)}~dx' ,
\label{Ladjoint}
\eeq
so that $\Ncal^\perp$ (rather than $\Ncal$) is invariant under $L^*$.

To summarise: $\Ncal$ is invariant under every element of $\Lcal$, and $\Ncal^\perp$
is invariant under every element of $\Rcal$ and under the set $\Lcal^*$ of adjoints.

{\bf Nest algebras.} The set of projections $P_a$ is a {\sl
nest}, that is, it is linearly ordered, contains $0$ and $1$ and is
order-complete (equivalently, closed in the strong operator topology). We use
the same notation $\Ncal$ for this set as for the underlying
subspaces, as it should not cause confusion.
The {\sl nest algebra} $\mathrm{Alg}\,\Ncal$ is the set of bounded operators
leaving invariant the elements of $\Ncal$. Thus $\mathrm{Alg}\,\Ncal$ contains the
bounded operators of type $L$.

The {\sl lattice} $\mathrm{lat}\,\Lcal$ of $\Lcal$ is the set of all subspaces
left invariant by every operator in $\Lcal$. A nest algebra $\Ecal$, say, has the form
$\Ecal = \mathrm{Alg}(\mathrm{lat}\,\Ecal)$. This property is called
{\sl reflexivity},
and it is shared with von Neumann algebras (strongly closed $*$-algebras of
operators).
Similarly a projection lattice is reflexive if it equals
$\mathrm{lat}$ of its algebra; nests are reflexive.

\begin{proposition}
\label{prop:closure}
The weak closure $\overline{\Lcal}$ of the bounded $L$-type operators is the nest
algebra $\mathrm{Alg}\,\Ncal$, that is, the set of all
operators leaving $\Ncal$ invariant. Consequently $\mathrm{lat}\,\Lcal=\Ncal$.
\end{proposition}

Throughout, the bar denotes the closure in the weak operator topology (where 
$T_\alpha \rightarrow T$ if $\langle T_\alpha f,g \rangle \rightarrow \langle T f,g \rangle$ 
for all $f, g \in H)$,
rather than the norm closure.

\begin{proof}
Rank-one elements of
$\mathrm{Alg}\,\Ncal$ have kernels $f(x)g(x')$ with $f$ supported in
$x>a$ and $g$ in
$x'<a$ for some $a$, so they are $L$-type. By Erd\H{o}s' theorem
\cite{erdos1968finite}, every finite-rank operator in a nest algebra is a sum of
rank-one operators in the algebra, and the finite-rank operators are weakly dense
in it. Hence $\overline{\Lcal}\supseteq\mathrm{Alg}\,\Ncal$, and the reverse
inclusion holds because $\mathrm{Alg}\,\Ncal$ is weakly closed and contains $\Lcal$.
Finally, $\mathrm{lat}\,\Lcal=\mathrm{lat}\,\overline{\Lcal} =\mathrm{lat}(\mathrm{Alg}\,\Ncal)=\Ncal$, since nests are reflexive.
\end{proof}

The set $\Lcal$ itself is not weakly closed. For example $\mathrm{Alg}\,\Ncal$
contains operators, such as translations $\psi(x)\mapsto\psi(x-1)$, whose kernels
are delta functions away from the diagonal. These are neither integral operators
with $L^2$ kernels nor multiplication operators, and so are not in $\Lcal$.

{\sl Finite-dimensional analogue.} The discretised analogues of the integral
operators in $\Lcal$ are the lower triangular matrices, which form a nest algebra.
With respect to the standard basis $\{e_i\}$, let $q_i$ be the projection onto
$e_i$. A single matrix entry $a_{ij}$ corresponds to a multiple of the rank-one
operator $c_{ij}=e_j\otimes e_i$ defined by
$(e_j\otimes e_i) b ~=~ \langle b,e_j\rangle e_i$, which has the
property $c_{ij} = c_{ij} q_j = q_i c_{ij}$. So the set of lower triangular matrices
consists of all linear combinations of such operators with $i\ge j$, and each entry
$a_{ij}$ corresponds to a discrete analogue of an integral operator with kernel
$K(x,x') = a_{ij}~\delta(x-x_i)~\delta(x'-x_j)$. Proposition~\ref{prop:closure} is
the continuous counterpart, with rank-one operators again playing the central role.

\medbreak
{\bf Inverses of $L$-type operators.}
An element of $\mathrm{Alg}\,\Ncal$ that is invertible in $B(H)$ need not have its
inverse in $\mathrm{Alg}\,\Ncal$. For example on $L^2(\mathbb{R})$, translation
$\psi(x)\mapsto\psi(x-1)$ is a unitary element of $\mathrm{Alg}\,\Ncal$ whose inverse
is not. For second-kind $L$-operators with Hilbert--Schmidt kernels (as holds for
standard scattering kernels on a bounded interval), however, the inverse is again of
$L$-type.

\begin{proposition}
\label{prop:inverse}
Let $L=f(I+V)$ on $H=L^2(\Omega)$, where $f$ is invertible in $L^\infty(\Omega)$ and
$V$ is a Volterra operator with a Hilbert--Schmidt kernel supported in $x'<x$ (for
example, on a bounded interval $\Omega=[0,\ell]$). Then $L$ is invertible,
\beq
L^{-1} ~=~ (I+V)^{-1}f^{-1} ~=~ \sum_{n\ge0}(-V)^n f^{-1},
\label{eq:Linverse}
\eeq
with the series converging in operator norm, and $L^{-1}\in\Lcal$, with diagonal
part $f^{-1}$.
\end{proposition}

\begin{proof}
Because $V$ is compact and leaves invariant the continuous nest $\Ncal$, it is
quasinilpotent (Ringrose's theorem; see \cite{gohberg1970volterra,davidson1988nest}),
so the Neumann series for $(I+V)^{-1}$ converges in operator norm. Each power $V^n$
with $n\ge1$ has a Hilbert--Schmidt kernel supported in $x'<x$, and
$\Vert{}V^n\Vert{}_{HS}\le\Vert{}V^{n-1}\Vert{}\,\Vert{}V\Vert{}_{HS}$ with $\Vert{}V^{n-1}\Vert{}^{1/n}\to0$, so the terms
with $n\ge1$ also converge in Hilbert--Schmidt norm. Their sum is an operator whose
kernel (the resolvent kernel of $V$) is supported in $x'<x$, so $L^{-1}$ is the
multiplication operator $f^{-1}$ plus an $L$-type integral operator.
\end{proof}

This is what allows forward marching:
 $Lu = f$ (each forward sweep of the L-R series) determines $u(x)$ sequentially from left to right, 
 using only values to the left of $x$.
 (The inverse of a lower triangular matrix
with non-zero diagonal is lower triangular for the same reason.) For typical
first-kind equations with smoothing kernels, such as logarithmic single-layer
kernels, the inverse is unbounded on $L^2$, since the operator is
compact; a formulation between appropriate Sobolev spaces may instead  
be well posed \cite{gohberg1970volterra}.

\medbreak
{\bf Triangularity.} In the sense of Kadison and Singer
\cite{kadison1960triangular}, a subalgebra $\Acal\subset B(H)$ is {\sl triangular} if
$\Acal\cap\Acal^*$ is maximal Abelian in $B(H)$, and {\sl maximal triangular} if it
is not properly contained in any larger triangular algebra. By (\ref{Ladjoint}),
$\Lcal\cap\Lcal^*$ consists precisely of the multiplication operators $L^\infty(\Omega)$,
which form a maximal Abelian algebra on $L^2(\Omega)$; hence $\Lcal$ is triangular.
The following is essentially classical
\cite{kadison1960triangular,ringrose1965some,davidson1988nest}; we include the short
argument:

\begin{proposition}
\label{prop:maxtri}
The weak closure $\overline{\Lcal}=\mathrm{Alg}\,\Ncal$ is maximal triangular.
\end{proposition}

\begin{proof}
An operator lies in $\mathrm{Alg}\,\Ncal\cap(\mathrm{Alg}\,\Ncal)^*$ exactly when it
commutes with every $P_a$, so this intersection is the commutant
$\Ncal'=(L^\infty)'=L^\infty$, and $\mathrm{Alg}\,\Ncal$ is triangular. Now let
$\Bcal\supseteq\mathrm{Alg}\,\Ncal$ be triangular, let $T\in\Bcal$, and put
$X=P_a^\perp T P_a$. Then $X\in\Bcal$, and $X^*=P_aT^*P_a^\perp$ maps functions
supported in $x\le a$ into $\Mcal_a$, so it leaves every $\Mcal_b$ invariant (for
$b\ge a$ it vanishes on $\Mcal_b$, and for $b<a$ its range lies in
$\Mcal_a\subseteq\Mcal_b$). Hence $X^*\in\mathrm{Alg}\,\Ncal\subseteq\Bcal$ and
$X\in\Bcal\cap\Bcal^*$. This algebra is Abelian and contains $P_a$, so
$X=P_a^\perp P_aX=0$. Since $a$ was arbitrary, $T\in\mathrm{Alg}\,\Ncal$.
\end{proof}

Nest algebras and maximal triangular algebras share further structural
properties: 
for example, bounded derivations of nest algebras are implemented by bounded
operators \cite{christensen1977derivations,spivack1985derivations},
and a similar 
result holds for maximal triangular algebras under suitable conditions
\cite{spivack2025remarks}. For two coupled surfaces (for example in a
rough 2D waveguide),  
however, the corresponding nest has
multiplicity $2$, its commutant $L^\infty\otimes M_2$ is non-Abelian, and the
resulting algebra is no longer triangular.

\section{Arveson's distance formula and the cut blocks}
\label{sec:arveson}

Recall that $P_a$ is the orthogonal projection onto $\Mcal_a$, the functions
vanishing for $x\le a$, so that $P_a^\perp=I-P_a$ projects onto the functions
supported in $x\le a$. For any nest, Arveson's distance formula
\cite{arveson1975interpolation} states that, for every $T\in B(H)$,
\beq
\mathrm{dist}(T,\mathrm{Alg}\,\Ncal) ~=~ \sup_a \|P_a^\perp T P_a\| .
\label{arveson}
\eeq
Here the distance is measured in operator norm,
$\mathrm{dist}(T,\mathrm{Alg}\,\Ncal)=\inf\{\|T-L'\|: L'\in\mathrm{Alg}\,\Ncal\}$,
so it is the error of the best approximation to $T$ by an operator that leaves every
subspace of the nest invariant.

We can think of each $a$ as a {\sl cut}, splitting the surface into the part $x\le a$
and the part $x>a$. The block $P_a^\perp T P_a$ maps functions supported in $x>a$ to
functions supported in $x\le a$. For an integral operator its kernel is $K(x,x')$
restricted to $x<a<x'$, so it describes how the field beyond the cut feeds back onto
the surface before it: the backward coupling across the cut.

One half of (\ref{arveson}) is elementary: For any $L'\in\mathrm{Alg}\,\Ncal$ we have
$P_a^\perp L' P_a=0$, so $P_a^\perp T P_a=P_a^\perp(T-L')P_a$ and hence
$\|P_a^\perp T P_a\|\le\|T-L'\|$. The main substance of Arveson's theorem is the reverse
inequality: some $L'\in\mathrm{Alg}\,\Ncal$ has error exactly equal to the
worst single-cut norm (the infimum is attained, since $\mathrm{Alg}\,\Ncal$ is weakly
closed, and bounded sets are weakly compact), so the corrections needed at different cuts do not accumulate.
In finite dimensions this reduces to Parrott's matrix-completion lemma \cite{parrott1978quotient} 
(see also Power \cite{power1980distance}), and Nehari's theorem \cite{nehari1957bounded} is the 
translation-invariant special case.

The formula is useful here because it is not tied to triangular
truncation. The L-R method uses 
one particular forward approximation, $L=\Pcal(A)$, with remainder
$R=\Qcal(A)$, and 
since $\Pcal(A)\in\mathrm{Alg}\,\Ncal$ we have
$\mathrm{dist}(A,\mathrm{Alg}\,\Ncal)\le\|R\|$. However, one could equally split
$A=L'+(A-L')$ with any $L'\in\mathrm{Alg}\,\Ncal$: for example a forward operator
arising from a one-way (parabolic) approximation, a truncation with a modified
diagonal, or an $L'$ chosen to make the remainder as small as possible. Arveson's
formula gives the exact distance to the full nest algebra of forward-only
operators, which is generally larger than the narrower class of $L$-type
integral operators considered here. Consequently, it gives a lower bound for
every forward-only approximation restricted to that narrower class: no such
approximation to $A$ can have a remainder smaller than the worst cut
block. Thus it provides a benchmark against which truncation, or any other
choice of $L'$, can be judged, and it shows how much could in
principle be gained by choosing $L'$ differently. It does not by
itself decide whether the resulting 
iteration converges, which depends on the spectral radius of $(A-L')L'^{-1}$; a
sufficient condition is $\|A-L'\|\,\|L'^{-1}\|<1$.

\medbreak
{\bf The kernel.} Take the hard (Neumann) surface $z=h(x)$ in 2D,
$G=\frac{i}{4}H_0^{(1)}(k|{\bf r}-{\bf r}'|)$, and the unnormalised normal
$(-h'(x'),1)$, so that $A=\frac12 I - K_D$, where $K_D$ is the double-layer operator
(the boundary operator obtained from the double-layer potential
$\int_\Gamma \partial_{n'}G({\bf r},{\bf r}')\,
\phi({\bf r}')\,dS'$, with $\partial_{n'}$ the source-point normal
derivative). The boundary limits of this potential contain a jump term 
$\pm \frac12 I$, where the sign depends on the normal and the side from which
the limit is taken. 
The unnormalised normal here incorporates the surface Jacobian.
The operator $K_D$ has kernel
\beq
K(x,x') ~=~ \frac{ik}{4}\,H_1^{(1)}(kR)\,\frac{D(x,x')}{R},\qquad
D(x,x') = h(x)-h(x')-h'(x')(x-x'),
\label{eq:dlkernel}
\eeq
and $R=\sqrt{(x-x')^2+(h(x)-h(x'))^2}$. (Sign conventions for the $\frac12$ and the
normal vary but they are irrelevant here since $I\in\mathrm{Alg}\,\Ncal$.)

\begin{itemize}
\item {\sl Near the diagonal:} $D\approx\frac12 h''(x')(x-x')^2$ and
$H_1^{(1)}(t)\approx -2i/(\pi t)$, so
$K(x,x')\to h''(x)/\big(4\pi(1+h'(x)^2)\big)$. The kernel is bounded for a $C^2$
surface, so the cut block has no singularity at the corner $x=x'=a$.
\item {\sl Far from the diagonal} ($s=\vert{}x-x'\vert{}$, $ks\gg1$), using
$\vert{}H_1^{(1)}(t)\vert{}\approx(2/\pi t)^{1/2}$,
\beq
|K(x,x')| ~\lesssim~ \Big(\frac{k}{8\pi}\Big)^{1/2}
\Big(2\|h\|_\infty\, s^{-3/2} + \|h'\|_\infty\, s^{-1/2}\Big).
\eeq
The height term is integrable; the slope term is not.
\end{itemize}

\begin{proposition}
\label{prop:schur}
Let $h\in C^2([0,\ell])$ and let $A=\frac12 I-K_D$ be the 2D Neumann boundary
integral operator with kernel (\ref{eq:dlkernel}). Then
\beq
\mathrm{dist}(A,\mathrm{Alg}\,\Ncal) ~\le~ \|R\| ~\le~
C\,\frac{\|h''\|_\infty}{k} + k\|h\|_\infty
+ \|h'\|_\infty\bigg(\Big(\frac{k\ell}{2\pi}\Big)^{1/2}+\frac{1}{20}\bigg),
\label{eq:schurbound}
\eeq
where $C$ depends only on $\|h'\|_\infty$.
\end{proposition}

\begin{proof}
Let $s=|x-x'|$. By the mean value theorem $R=s\,(1+h'(\xi)^2)^{1/2}$ for some
$\xi$, so $s\le R\le\tau s$ with $\tau=(1+\|h'\|_\infty^2)^{1/2}$. Since the remainder operator
$R=\Qcal(A)$ is strictly upper triangular, Schur's test bounds both $\|R\|$ and every cut block
$\|P_a^\perp AP_a\|$ by $\int_0^\ell\kappa(s)\,ds$, where $\kappa(s)$ is the supremum of
$|K(x,x')|$ over $|x-x'|=s$; and $\mathrm{dist}(A,\mathrm{Alg}\,\Ncal)\le\|R\|$ by
(\ref{arveson}).

{\sl Near zone, $ks\le1$:} Here $0<kR\le\tau$, and $|H_1^{(1)}(t)|\le C_\tau/t$ on
$(0,\tau]$. With $|D|\le\frac12\|h''\|_\infty s^2\le\frac12\|h''\|_\infty R^2$ this gives
$|K|\le\frac18C_\tau\|h''\|_\infty$, so this zone contributes at most
$C\|h''\|_\infty/k$.

{\sl Far zone, $ks>1$:} Here $kR\ge ks>1$. For $t\ge1$,
$|H_1^{(1)}(t)|\le(2/\pi t)^{1/2}\big(1+1/(5t^2)\big)$ (the ratio
$(\pi t/2)^{1/2}|H_1^{(1)}(t)|$ equals $1.124$ at $t=1$ and tends to $1$ like
$1+3/(16t^2)$; the bound is readily verified numerically). With $|D|/R\le(2\|h\|_\infty+\|h'\|_\infty s)/s$ this gives
\[
|K| ~\le~ \Big(\frac{k}{8\pi}\Big)^{1/2}\Big(1+\frac{1}{5k^2s^2}\Big)
\Big(2\|h\|_\infty s^{-3/2}+\|h'\|_\infty s^{-1/2}\Big).
\]
Integrating over $1/k<s<\ell$, using $\int s^{-1/2}<2\ell^{1/2}$,
$\int s^{-3/2}<2k^{1/2}$, $\int s^{-5/2}<\frac23k^{3/2}$ and
$\int s^{-7/2}<\frac25k^{5/2}$, the slope terms contribute less than
$\|h'\|_\infty\big((k\ell/2\pi)^{1/2}+0.03\big)$ and the height terms less than
$0.84\,k\|h\|_\infty$. Adding the two zones gives (\ref{eq:schurbound}).
\end{proof}

These bounds have two consequences. First, the Schur test gives the same upper
bound for the cut blocks as for $R$, so it cannot distinguish the Arveson distance
from the truncation remainder.
Membership of the Schur class alone does not guarantee a uniform comparison $\Vert{}R\Vert{}\le C\,\mathrm{dist}(A,\mathrm{Alg}\,\Ncal)$ in $L^2$ operator norm: piecewise-constant kernels 
on $[0,\ell]$ belong to the Schur class and represent $n\times n$ matrices,
for which $\Vert{}R\Vert{}$ can exceed $\Vert{}A\Vert{}$ (and hence 
$\mathrm{dist}(A,\mathrm{Alg}\,\Ncal)$) by a factor of order $\log n$~\cite{kwapien1970main}.
 A comparison with constants depending on additional restrictions
on physical Green's-function kernels remains an open question
(Section~\ref{sec:questions}).
Second, the $(k\ell)^{1/2}$ growth comes from the $s^{-1/2}$ tail of the slope
term. This is consistent with the known difficulties of the free-space Green's
function on infinite rough surfaces.

\medbreak
{\bf Numerical check.} The operators were discretised by the Nystr\"om
method on a uniform grid with spacing $\Delta x=\lambda/10$ ($\lambda=1$).
Operator norms are matrix $2$-norms, and the reported distances are maxima
over 25 sampled cut positions; thus they are lower bounds for the discrete
Arveson distances obtained by maximising over all cuts. Random surfaces were
generated by spectral synthesis with fixed random seeds, using Gaussian spectra
for the length and slope experiments and a power-law spectrum for the comparison
below. Each result is for a single realisation, not an ensemble average. In the
length-dependence experiment a new realisation is generated at each length,
rather than extending one fixed surface. The nominal rms height is set on a
longer periodic generating surface before the required segment is extracted.
For the following table the nominal rms height is $0.25\lambda$, the Gaussian
correlation length is $2\lambda$, and the measured rms slopes are approximately
$0.16$--$0.23$:
\begin{center}
\begin{tabular}{lcccc}
 & dist & $\|R\|$ & $\|B\|$ & $\rho(B)$\\ \hline
$\ell=10\lambda$ & 0.32 & 0.35 & 0.70 & 0.04\\
$\ell=40\lambda$ & 0.43 & 0.58 & 1.15 & 0.17\\
$\ell=160\lambda$ & 1.11 & 1.49 & 2.98 & 0.34
\end{tabular}
\end{center}
\begin{itemize}
\item For these three surface realisations, the distance is consistent
with growth of order $\ell^{1/2}$ and is $0.7$--$0.9$ of $\Vert{}R\Vert{}$;
truncation is relatively close to the best forward approximation
in these examples.
\item At fixed $\ell=40\lambda$ the ratio of distance to rms slope is
approximately $2.7$ across the amplitudes tested. This is consistent with
dominance of the slope-dependent contribution over this range, for the
particular surface realisation used.
\item The spectral radius $\rho(B)$ is much smaller than $\Vert{}B\Vert{}$, so norm estimates
badly overstate the convergence problem.
\end{itemize}

{\bf The extremal input.} For the Gaussian surface with $\ell=80\lambda$ (rms
height $0.25\lambda$, correlation length $2\lambda$), the cut block at the
mid-point is approximately rank one (largest singular value $0.69$, compared with $0.08$ for the next), 
so its norm is attained to good 
approximation by a single input function: its right singular vector. This is a
predominantly left-going grazing wave ($k_x\approx-k$) spread over tens of wavelengths
from the cut. 
Applied to unit-norm windowed right-going waves $\psi(x)\propto e^{ikx\cos\theta}$, 
the output norm $\Vert{}P_a^\perp A P_a \psi\Vert{}$ is $\sim 10^{-4}$ for 
$0^\circ\le\theta\le30^\circ$ on the Gaussian surface (negligible spectral 
content at the Bragg wavenumber $2k$), against $0.02$--$0.2$ for left-going waves. 
On a power-law surface with Bragg content (rms height $0.1\lambda$) the corresponding 
norms are $\sim 0.01$ for right-going and $0.01$--$0.06$ for left-going waves.
So on all of $L^2$ the Arveson distance mostly measures a left-going wave
propagating across the cut, which is forward scattering for that wave. The
coupling relevant to predominantly right-going incident fields can therefore be
much smaller than the unrestricted Arveson distance suggests. This motivates
considering the cut block restricted to right-going fields, which measures
Bragg backscatter, and seeking a directional analogue of the distance formula.

\section{The mirror surface}
\label{sec:mirror}

By Proposition~\ref{prop:inverse}, for a given surface $L^{-1}$ is itself an
$L$-type integral operator in the second-kind case. It is natural to ask whether it
is the $L$-operator of some equivalent `virtual surface'
(Question~\ref{q:virtual} in Section~\ref{sec:questions}).
The exact form of this question     
probably has a negative answer: a boundary-integral
kernel $K_h(x,x')$ depends on the surface only through $h$ and $h'$ at the two
points $x,x'$, whereas the resolvent kernel in $L^{-1}$ depends on $h$ along the
whole interval between them. For the hard surface, however, there is an exact
identity with a natural virtual surface, namely the mirror image $z=-h(x)$.

\begin{proposition}[Mirror surface]
\label{prop:mirror}
Let $L_h=\Pcal(A_h)=\frac12 I-K^L_h$ be the
$L$-operator of the hard surface $z=h(x)$ on a bounded interval, where $K^L_h$ is
the lower part of the double-layer operator with kernel
(\ref{eq:dlkernel}). Then
\beq
K^L_{-h}~=~-K^L_h, \qquad
L_h L_{-h} ~=~ L_{-h} L_h ~=~ \tfrac14 I-(K^L_h)^2 ,
\label{eq:mirror1}
\eeq
and
\beq
L_h^{-1} ~=~ 4\,L_{-h}\big(I-4(K^L_h)^2\big)^{-1}
~=~ 4\,L_{-h}\sum_{n\ge0}\big(4(K^L_h)^2\big)^n .
\label{eq:mirror2}
\eeq
\end{proposition}

\begin{proof}
In the kernel, $D(x,x')=h(x)-h(x')-h'(x')(x-x')$ is linear in $h$,
while $R$, and hence $H_1^{(1)}(kR)/R$, depends on $h$ only through
$(h(x)-h(x'))^2$. The diagonal value $h''/\big(4\pi(1+h'^2)\big)$ is also odd in $h$.
So the kernel is odd in $h$ pointwise, and since $\Pcal$ is linear,
$K^L_{-h}=-K^L_h$. Then $L_{-h}=\frac12 I+K^L_h$, which commutes with
$L_h=\frac12 I-K^L_h$, and multiplying out gives (\ref{eq:mirror1}). The operator
$K^L_h$ is a Volterra operator, and thus quasinilpotent, so $4(K^L_h)^2$ is also
quasinilpotent. Hence $I-4(K^L_h)^2$ is always invertible and the series in
(\ref{eq:mirror2}) converges in operator norm, even when $\|K^L_h\|\ge\frac12$.
\end{proof}

{\bf Remarks.}
\begin{itemize}
\item The correction is exactly
$L_h^{-1}-4L_{-h}=16\,L_{-h}(K^L_h)^2\big(I-4(K^L_h)^2\big)^{-1}$, which begins at second order in $K^L_h$ (the factor
$L_{-h}=\frac12 I+K_h^L$ also contributes odd powers). Since the kernel is small for gentle
slopes, to first order in slope the inverse of the forward operator is (four times)
the forward operator of the reflected surface.
\item The identity is equivalent to splitting the Neumann series for $L_h^{-1}$ into
its even and odd parts:
\beq
\begin{split}
\big(\tfrac12 I-K^L_h\big)^{-1} &= 2\sum_{n\ge0}\big(2K^L_h\big)^n
~=~ 2\big(I+2K^L_h\big)\sum_{n\ge0}\big(4(K^L_h)^2\big)^n \\
&= 4\,L_{-h}\sum_{n\ge0}\big(4(K^L_h)^2\big)^n .
\end{split}
\eeq
The odd-order terms are thus generated by the forward operator of the mirror
surface, and replacing $L_h^{-1}$ by $4L_{-h}$ amounts to keeping the first two
terms of the series, i.e.\ the first-order Born approximation to the forward solution.
The significance of the proposition is therefore this physical interpretation, made
possible by the parity $K_{-h}=-K_h$, rather than a new computational device.
\item Computationally, $L_{-h}=\frac12 I+K^L_h$ is constructed from the same matrix
$K^L_h$ as $L_h$, so no new surface calculation is needed. In 2D, applying
$L_h^{-1}$ to a vector is a forward substitution (the forward sweep of the L-R
method), at the same $O(N^2)$ cost as a matrix--vector product, and spectral
acceleration reduces both further. Replacing $L_h^{-1}$ by $4L_{-h}$ in the series
therefore saves nothing asymptotically in 2D, and it introduces an error of second
order in slope. The possible gain in 3D is discussed below.
\item A numerical check with the discretisation of Section~\ref{sec:arveson}
(Gaussian-spectrum surfaces, $\ell=40\lambda$, correlation length $2\lambda$) gives
\begin{center}
\begin{tabular}{ccc}
rms slope & $\|K^L_h\|$ & $\|L_h^{-1}-4L_{-h}\|/\|L_h^{-1}\|$\\ \hline
0.03 & 0.13 & 0.013\\
0.06 & 0.26 & 0.045\\
0.16 & 0.64 & 0.22
\end{tabular}
\end{center}
with (\ref{eq:mirror1}) holding to machine precision. The relative error grows
roughly quadratically in slope.
\item The parity argument applies unchanged in 3D, where
$D=h({\bf x})-h({\bf x}')-\nabla h({\bf x}')\cdot({\bf x}-{\bf x}')$ is again
linear in $h$, so $K_{-h}=-K_h$ and (\ref{eq:mirror1}) still hold; nor does it
depend on the sign convention $A=\frac12 I\pm K$. If the surface is ordered by $x$,
however, a discretisation into finite-width $x$-slices produces a
block lower-triangular matrix with nontrivial diagonal blocks, coupling
the transverse degrees of freedom within each slice. This does not by itself imply that the continuous
$x$-ordered integral operator loses its Volterra or quasinilpotent
character. 
A sufficient condition is that the kernel norm be bounded by an integrable function 
of $x-x'$ alone: if $\Vert{}K(x,x')\Vert{}\le k(x-x')$ for $x'<x$, where the norm is 
that of operators in the transverse variable and $k\in L^1(0,\ell)$, then
by Young's inequality
$\|V^n\|\le\|k^{*n}\|_{L^1(0,\ell)}$, and
$\|k^{*n}\|_{L^1(0,\ell)}\le e^{\sigma\ell}\big(\int_0^\ell e^{-\sigma t}|k(t)|\,dt\big)^n$
for every $\sigma>0$. Letting $\sigma\to\infty$ shows that $V$ is quasinilpotent.
For the weakly singular 3D double layer one might expect $k(s)\sim\log(1/s)$, but this,
and the precise function spaces involved, still need to be verified. In particular, mere integrability of an
operator-valued kernel does not justify the factorial bound
$\|V^n\|\le C^n/n!$ stated in some elementary Volterra estimates.
If quasinilpotence holds, then (\ref{eq:mirror2}) applies in the
continuous setting as well. The nontrivial diagonal blocks in a
finite-width slice discretisation must be analysed separately.
For a discretised operator, the
diagonal blocks can affect the spectrum, and convergence of the series
in (\ref{eq:mirror2}) must be checked, for example through
$\rho(2K_h^L)<1$. Where that condition holds, replacing $L_h^{-1}$
by $4L_{-h}$ plus a few correction terms might avoid expensive
diagonal block solves. Whether this is advantageous depends on the
discretisation and requires numerical testing.
\item The identity uses the second-kind (Neumann) formulation with a constant
diagonal term; the Dirichlet (single-layer, first-kind) case needs separate
treatment.
\end{itemize}

\section{Further questions}
\label{sec:questions}

\begin{enumerate}
\item \label{q:virtual}
{\sl Virtual surfaces.} In the second-kind setting, $L^{-1}$ is again an
$L$-type operator (Proposition~\ref{prop:inverse}). Under what conditions can
it be realised as the forward operator associated with a `virtual surface'?
More specifically, can the correction
$\big(I-4(K^L_h)^2\big)^{-1}$ in Proposition~\ref{prop:mirror}
be incorporated into a modified surface, Green's function or boundary
condition? This may be related to inverse marching methods \cite{SR1,SR2,chen2018rough}.

\item {\sl Optimality of truncation.} For physically restricted classes of
rough-surface Green's-function kernels, can one bound the norm of the
backward triangular part by
$\|R\|\le C\,\mathrm{dist}(A,\mathrm{Alg}\,\Ncal)$,
with $C$ controlled by geometric or spectral parameters? Schur-class
membership alone is insufficient, but the observed near-rank-one structure
of the worst cut block may help identify suitable restrictions.

\item {\sl A directional distance formula.} The position nest
$\{\Mcal_a\}$ and the wavenumber nest of shifted Hardy spaces
$e^{i\lambda x}H^2(\mathbb{R})$ form the Fourier binest of Katavolos
and Power \cite{katavolos1997fourier}, related to the parabolic algebra
\cite{katavolos2002translation,kastis2017parabolic}.
Does its algebra admit an Arveson-type distance formula, and could such
a formula quantify the restricted cut-block norms
$\sup_a\|P_a^\perp A P_a Q_0\|$, where $Q_0$ projects onto
right-going wavenumbers? Distance formulae do not hold for all
subspace lattices.

\item {\sl Three dimensions.} In three dimensions, finite-width
$x$-slices lead to nontrivial diagonal blocks even when the continuous
$x$-ordered operator is Volterra. Under what conditions do the
continuous nest-algebra properties and the mirror identity lead to
stable, efficient discrete marching schemes (Section~\ref{sec:mirror})?

\item {\sl First-kind equations and unbounded operators.} For the
Dirichlet single-layer formulation, the inverse need not be bounded on
$L^2$. Can the results be formulated on appropriate Sobolev spaces,
and is there a counterpart of the mirror identity?
\end{enumerate}

\bibliography{references}

\end{document}